\input amstex
\documentstyle{amsppt}
\input bull-ppt
\define\shiftsum#1#2#3#4{\sum\limits_{#1\,_#4\le#2_#4%
\le#3_#4\atop#4=1,2,
\dots,\ell}}
\define\qrfac#1#2{\left(#1\right)_{#2}} 
\define\smprod{\prod\limits^{\ell}_{k=1}}
\define\multsum#1#2#3{\sum\limits_{0\le #1_#3 \le #2_#3 
\atop #3 =1,2,\dots,\ell}}
\define\bcl{\bmi{C_\ell}}
\define\bal{\bmi{A_\ell}}
\define\sqprod{\prod\limits_{r\!,s=1}^{\ell}}
\define\xover#1#2{{x_#1\over x_#2}}

\topmatter
\cvol{26}
\cvolyear{1992}
\cmonth{April}
\cyear{1992}
\cvolno{2}
\cpgs{258-263}
\title The $\bal$\  and $\bcl$\  Bailey Transform and 
Lemma\endtitle
\shorttitle{THE $A_\ell$ AND $C_\ell$ BAILEY TRANSFORM AND 
LEMMA}
\author Stephen C. Milne and Glenn M. Lilly\endauthor
\shortauthor{S. C. Milne and G. M. Lilly}
\address Department of Mathematics, 
The Ohio State University, Columbus, Ohio 43210\endaddress
\ml milne\@function.mps.ohio-state.edu \endml
\thanks S. C. Milne was partially supported by NSF grants 
DMS 86-04232, 
DMS 89-04455, and DMS 90-96254\endthanks
\thanks G. M. Lilly was fully supported by NSA supplements 
to the above 
NSF grants and by NSA grant MDA 904-88-H-2010\endthanks
\date April 28, 1991\enddate
\abstract 
We announce a higher-dimensional generalization of the 
Bailey
Transform, 
Bailey Lemma, and iterative ``Bailey chain'' concept 
in the setting of basic hypergeometric series very 
well-poised
on unitary $A_{\ell}$ or symplectic $C_{\ell}$ groups.  
The classical case, corresponding to $A_1$ or equivalently 
$\roman U(2)$, 
 contains an immense amount of the theory and application of
one-variable basic hypergeometric series, including 
elegant proofs
 of the Rogers-Ramanujan-Schur identities.  In particular, 
our 
program extends much of the classical work of Rogers, 
Bailey, 
Slater, Andrews, and Bressoud. \endabstract
\endtopmatter

\document

\heading 1. Introduction\endheading

The purpose of this paper is to announce a 
higher-dimensional generalization of
the Bailey Transform [2] and Bailey Lemma [2] in the 
setting of basic
hypergeometric series very well-poised on unitary [19] or 
symplectic [14]
groups.  Both types of series are directly related [14, 
18] to the
corresponding Macdonald identities. The series in [19] 
were strongly motivated
by certain applications of mathematical physics and the 
unitary groups 
$\roman U(n)$
in [10, 11, 15, 16].  The unitary series use the notation 
$A_\ell$, or
equivalently $\roman U
(\ell +1)$; the symplectic case, $C_\ell$.  The classical 
Bailey
Transform, Lemma, and very well-poised basic 
hypergeometric series correspond
to the case $A_1$, or equivalently $\roman U(2)$.

The classical Bailey Transform and Bailey Lemma contain an 
immense
amount of the theory and application of one-variable basic 
hypergeometric
series [2, 12, 25].  They were ultimately inspired by 
Rogers' [24] second proof
of the Rogers-Ramanujan-Schur identities [23].  The Bailey 
Transform was first
formulated by Bailey [8], utilized by Slater in [25], and 
then recast by
Andrews [4] as a fundamental matrix inversion result.  
This last version of the
Bailey Transform has immediate applications to connection 
coefficient theory
and ``dual'' pairs of identities [4], and $q$-Lagrange 
inversion and quadratic
transformations [13].

The most important application of the Bailey Transform is 
the Bailey Lemma.  This result was
mentioned by Bailey [8; \S 4], and he described how the 
proof would work.  However, he
never wrote the result down explicitly and thus missed the 
full power of {\it iterating} it. 
Andrews first established the Bailey Lemma explicitly in 
[5] and realized its numerous possible
applications in terms of the iterative ``Bailey chain'' 
concept.  This iteration mechanism
enabled him to derive many $q$-series identities by 
``reducing'' them to more
elementary ones.  For example, the Rogers-Ramanujan-Schur 
identities can be reduced to
the $q$-binomial theorem.  Furthermore, general multiple 
series Rogers-Ramanujan-Schur
identities are a direct consequence of iterating suitable 
special cases of Bailey's Lemma. 
In addition, Andrews notes that Watson's $q$-analog of 
Whipple's transformation is an
immediate consequence of the second iteration of one of 
the simplest cases of Bailey's
Lemma.  Continued iteration of this same case yields 
Andrews' [3] infinite family of
extensions of Watson's $q$-Whipple transformation. Even 
Whipple's original work [26, 27]
fits into the $q=1$ case of this analysis.  Paule [22] 
independently discovered important
special cases of Bailey's Lemma and how they could be 
iterated.  Essentially all the depth of
the Rogers-Ramanujan-Schur identities and their iterations 
is embedded in Bailey's Lemma.

The process of iterating Bailey's Lemma has led to a wide 
range of applications
in additive number theory,  combinatorics, special 
functions, and mathematical
physics.  For example, see [2, 5, 6, 7, 9].

The Bailey Transform is a consequence of the terminating 
${}_4\phi_3$ summation
theorem.  The Bailey Lemma is derived in [1] directly from 
the ${}_6\phi_5$
summation and the matrix inversion formulation  [4, 13] of 
the Bailey
Transform. We employ a similar method in the $A_\ell$ and 
$C_\ell$ cases by
starting with a suitable, higher-dimensional, terminating 
${}_6\phi_5$
summation theorem extracted from [19] and [14], 
respectively.  The $A_\ell$
proofs appear in [20, 21], and the $C_\ell$ \ case is 
established in [17]. 
Many other consequences of the $A_\ell$ and $C_\ell$ 
generalizations of
Bailey's Transform and Lemma will appear in future papers. 
 These include
$A_\ell$ and $C_\ell$\ $q$-Pfaff-Saalsch\"utz summation 
theorems, $q$-Whipple
transformations, connection coefficient results, and 
applications of iterating
the $A_\ell$ or $C_\ell$ Bailey Lemma.

 \heading 2. Results\endheading
Throughout this article, let $\bmi{i}$, $\bmi{j}$, $\bmi{ 
N}$,
 and $\bmi{y}$ be vectors of length $\ell$ 
with nonnegative integer components.  Let $q$ be a complex 
number such that $|q|<1$.  Define
$$
\qrfac\alpha\infty\equiv(\alpha;\  q)_\infty:=
\prod\limits_{k\ge 0}(1-\alpha q^k)\tag2.1a$$
and, thus,
$$
\qrfac\alpha{n}\equiv(\alpha;\  q)_n:
=\qrfac\alpha\infty/\qrfac{\alpha q^n}\infty.\tag2.1b$$

Define the Bailey transform matricies, $M$ and $M^\ast$, 
as follows.

\dfn\nofrills{Definition \RM{($\bmi{M}$ and $\bmi{M}^\ast$ 
for $\bal$).} }
\ Let $a, x_1,\ldots,x_\ell$ be indeterminate.
Suppose that none of the denominators in (2.2a--b) 
vanishes. Then let
$$M(\bmi{ i}; \bmi{j};
A_\ell):= \sqprod\qrfac{q\xover rs 
q^{j_r-j_s}}{i_r-j_r}^{-1}\quad
\smprod\qrfac{aq\xover k\ell}{i_k+(j_1+\cdots+
j_\ell)}^{-1}; \tag2.2a$$
and
$$\align
& M^\ast(\bmi{ i};\ \bmi{j};\  A_\ell)\tag2.2b\cr
&\qquad \coloneq\smprod \biggl[1-a\xover k\ell
q^{i_k+(i_1+\cdots+i_\ell)}\biggr] \smprod 
\qrfac{aq\xover k\ell}{j_k+(i_1+\cdots+i_\ell)-1}\cr
&\quad\qquad\times \sqprod\qrfac{q\xover rs q^{j_r-j_s}}
{i_r-j_r}^{-1}\ (-1)^{(i_1+\cdots+i_\ell)-(j_1+\cdots+
j_\ell)}
\ q^{(i_1+\cdots+i_\ell)-(j_1+\cdots+j_\ell)\choose 2}.\cr
\endalign$$
\enddfn

\dfn\nofrills{Definition {\rm($\bmi{M}$ and $\bmi{ 
M}^\ast$ for $\bcl$).} }
\ Let $x_1,\ldots,x_\ell$ be indeterminate.
Suppose that none of the denominators in (2.3a--b) 
vanishes. Then let
$$
M(\bmi{i}; \bmi{j};  C_\ell):=
\sqprod\Biggl[\qrfac{q\xover rs 
q^{j_r-j_s}}{i_r-j_r}^{-1}\qrfac{qx_rx_sq^{j_r+
j_s}}{i_r-j_r}^{-1}\Biggr];
\tag2.3a
$$
and
$$
\align
\noalign{\hbox{(2.3b)}}
&M^\ast(\bmi{ i}; \bmi{ j};
  C_\ell)\cr
&\qquad \coloneq\sqprod\Biggl[\qrfac{q\xover rs 
q^{j_r-j_s}}{i_r-j_r}^{-1}
\qrfac{x_rx_sq^{j_r+i_s}}{i_r-j_r}^{-1}\Biggr]
\prod_{1\le r<s\le\ell}\left[{1-x_rx_sq^{j_r+j_s}
\over 1-x_rx_sq^{i_r+i_s}}\right]\cr 
&\qquad\quad\ \times (-1)^{(i_1+\cdots+i_\ell)-
(j_1+\cdots+j_\ell)} q^{(i_1+\cdots+i_\ell)-(j_1+\cdots+
j_\ell)\choose
2}\. \endalign$$
\enddfn

As in the classical case [1], we have the following theorem.
\proclaim\nofrills{Theorem {\rm(Bailey Transform for 
$\bal$ \ and $\bcl$)}\RM.}
\ Let $G=A_\ell$ or $C_\ell$.  Let $M$ and $M^\ast$ be 
defined
as in {\rm(2.2)} and {\rm(2.3)}, with rows and columns 
ordered lexicographically.  
Then $M$ and $M^\ast$ are inverse, infinite, 
lower-triangular matricies.  
That is,
$$\smprod\delta(i_k, j_k)=\shiftsum jyik M(\bmi{i};
\ \bmi{y};\  G)
\  M^\ast(\bmi{y};\ \bmi{j};\  G),\tag2.4$$
where $\delta(r,  s)=1$ if $r=s$, and $0$ otherwise.
\endproclaim

Equations (2.2) and (2.3) motivate the definition of the 
$A_\ell$ and $C_\ell$ Bailey pair.

\dfn\nofrills{Definition 
{\rm($\bmi{ G}$-Bailey Pair).}}\  
Let $G=A_\ell$ or $C_\ell$.  Let $N_k\ge0$ be integers for
$k=1,2,\ldots,\ell$.  Let $A=\{A_{(\bmi{ y};  G)}\}$ and 
$B=\{B_{(\bmi{y};
\  G)}\}$ be sequences. Let $M$ and $M^\ast$ be as above.
Then we say that $A$ and $B$ form a $G$-Bailey Pair
if
$$B_{(\bmi{ N};  G)}=\multsum yNk M(\bmi{ N};
\bmi{ y};  G)\  
A_{(\bmi{ y};  G)}.\tag2.5$$
\enddfn
As a consequence of the Bailey transform, 
(2.4), and the definition of the $G$-Bailey pair, (2.5),
we have the following result.

\proclaim\nofrills{Corollary \RM{(Bailey Pair 
Inversion).}}\ $A$ and $B$
 satisfy equation \RM{(2.5)} if and only if
$$A_{(\bmi{ N};  G)}=\multsum yNk M^\ast(\bmi{ N};
 \bmi{y};  G)\  B_{(\bmi{y};  G)}.\tag2.6$$
\endproclaim
Define the sequences $A'=\{A'_{(\bmi{ y};
  A_\ell)}\}$ and $B'=\{B'_{(\bmi{ y};  A_\ell)}\}$
by
$$\aligned
A'_{(\bmi{N};  A_\ell)}
:=& \smprod\qrfac{{aq\over\rho}\xover k\ell}{N_k}^{-1}
\quad\smprod\qrfac{\sigma\xover k\ell}{N_k}\cr
&\times
{\qrfac\rho{N_1+\cdots+N_\ell}\over\qrfac{aq/\sigma}{N_1+
\cdots+N_\ell}}
 (aq/\rho\sigma)^{N_1+\cdots+N_\ell}\ 
A_{(\bmi{ N};\  A_\ell)}\endaligned\tag2.7a
$$
and
$$\align 
 B'_{(\bmi{N};  A_\ell)}
:=\multsum yNk&\Biggl\{\smprod\left[\qrfac{\sigma\xover 
k\ell}{y_k}
\qrfac{{aq\over\rho}\xover k\ell}{N_k}^{-1}\right]
\sqprod\qrfac{q\xover rs 
q^{y_r-y_s}}{N_r-y_r}^{-1}\tag2.7b\\
&\qquad
\times{\qrfac{aq/\rho\sigma}{(N_1+\cdots+N_\ell)-(y_1+
\cdots+y_\ell)}\ 
\qrfac\rho{y_1+\cdots+y_\ell}\over
\qrfac{aq/\sigma}{N_1+\cdots+N_\ell}}\\
\shoveright{\times (aq/\rho\sigma)^{y_1+\cdots+y_\ell}\ B_{(
\bmi y;\  A_\ell)}\Biggr\}}
\endalign$$
Define the sequences $A'=\{A'_{(\bmi{y};
  C_\ell)}\}$ and $B'=\{B'_{(\bmi{y};  C_\ell)}\}$
by
$$ A'_{(\bmi{ N};  C_\ell)}:
=\smprod\left[{\qrfac{\alpha x_k}{N_k}\ 
\qrfac{qx_k\beta^{-1}}{N_k}\over
\qrfac{\beta x_k}{N_k}\ \qrfac{qx_k\alpha^{-1}}{N_k}}\right]
\left({\beta\over\alpha}\right)^{N_1+\cdots+N_\ell}
A_{(\bmi{N};  C_\ell)}\tag2.8a$$
and
$$\align
\noalign{\hbox{(2.8b)}}
B'_{(\bmi N;\  C_\ell)}:=
\multsum yNk&\Biggl\{cr\smprod\left[{\qrfac{\alpha 
x_k}{y_k}\
\qrfac{qx_k\beta^{-1}}{y_k}\over \qrfac{\beta x_k}{N_k}\ 
\qrfac{qx_k\alpha^{-1}}{N_k}}\right]
 \sqprod\qrfac{q\xover rs
q^{y_r-y_s}}{N_r-y_r}^{-1}\cr 
&\qquad\times\prod\limits_{1\le r<s\le\ell}
\left[\qrfac{qx_rx_sq^{y_r+y_s}}{N_s-y_s}^{-1}
\qrfac{qx_rx_sq^{N_s-y_s}}{N_r-y_r}^{-1}\right]\!{}\cr
&\qquad\times\qrfac{\beta\over\alpha}{(N_1+\cdots+N_\ell)
-(y_1+\cdots+y_\ell)}\ 
\left({\beta\over\alpha}\right)^{y_1+\cdots+y_\ell}
B_{(\bmi y;\  C_\ell)}\Biggr\}\cr
\endalign$$
These definitions lead to our generalization  of Bailey's 
lemma.

\proclaim\nofrills{Theorem 
\RM{(The $G$\<-generalization of Bailey's Lemma).}}\ Let 
$G=A_\ell$ or $C_\ell$.  Suppose $A=\{A_{(\bmi N;  G)}\}$
and $B=\{B_{(\bmi N;  G)}\}$ form a $G$-Bailey Pair.  If 
$A'=\{A'_{(\bmi N;
  G)}\}$
and $B'=\{B'_{(\bmi N; G)}\}$ are as above, 
then $A'$ and $B'$ also form a $G$-Bailey Pair.\endproclaim

\heading3. Sketches of Proofs\endheading

\demo{Proof of \RM{(2.4)}} 
In each case, $A_\ell$ and $C_\ell$, we begin with a 
terminating ${}_4\phi_3$ summation theorem. In the 
$C_\ell$ case, it
is first necessary to specialize Gustafson's $C_\ell\ 
_6\psi_6$ summation theorem, see [14], terminate it from 
below and then from above, and
further specialize the resulting terminating ${}_6\phi_5$ 
to yield a terminating ${}_4\phi_3$.  In both the $A_\ell$ 
and $C_\ell$ cases, the
${}_4\phi_3$ is modified by multiplying both the sum and 
product sides by some additional factors.  Finally, that 
result is transformed
term-by-term to yield the sum side of (2.4).\qed\enddemo

\demo{Proof of \RM{(2.6)}} Equation (2.6) follows directly 
from the definition, (2.5), and the termwise nature of the 
calculations in the proof of
(2.4).\qed\enddemo

\demo{Proof of Bailey's Lemma} The definitions in (2.7) 
and (2.8) are substituted into (2.5).  After an 
interchange of summation, the 
inner sum is
seen to be a special case of the appropriate ${}_6\phi_5$. 
 The ${}_6\phi_5$ is then summed, and the desired result
follows.\qed\enddemo 

Detailed proofs of the $C_\ell$ case will appear in [17], 
as will a discussion of the $C_\ell$ Bailey
chain and a connection coefficient result associated with 
the $C_\ell$ Bailey Transform.

\enddocument